\documentclass[12pt]{amsart}
\usepackage{amsmath,amstext,amsfonts,amsthm,amssymb}
\usepackage{eucal}
\usepackage{graphicx}
\renewcommand{\subsubsection}[1]{\addtocounter{subsubsection}{1}
{\ \\[3pt]\bf \thesubsubsection. \  #1} }

\swapnumbers

{  \theoremstyle{definition}

}
\newcommand{\Dist}{\operatorname{Dist}}

\newcommand{\Fil}{\operatorname{Fil}}

\newcommand{\Spec}{\operatorname{Spec}}

\newcommand{\hra}{\hookrightarrow}
\newcommand{\iso}{\overset{\sim}{\longrightarrow}}
\newcommand{\isom}{\overset{\sim}{=}}
\newcommand{\lra}{\longrightarrow}

\newcommand{\bea}{\begin{eqnarray*}}
\newcommand{\eea}{\end{eqnarray*}}
\newcommand{\bean}{\begin{eqnarray}}
\newcommand{\eean}{\end{eqnarray}}

\newcommand{\fsl}{\mathfrak{sl}}

\newcommand{\CR}{\mathcal{R}}

\newcommand{\BA}{\mathbb{A}}
\newcommand{\BC}{\mathbb{C}}

\newcommand{\BF}{\mathbb{F}}
\newcommand{\BG}{\mathbb{G}}

\newcommand{\BL}{\mathbb{L}}

\newcommand{\BQ}{\mathbb{Q}}
\newcommand{\BR}{\mathbb{R}}
\newcommand{\BS}{\mathbb{S}}
\newcommand{\BZ}{\mathbb{Z}}

\newcommand{\bfH}{\mathbf{H}}
\newcommand{\bfS}{\mathbf{S}}

\newcommand{\nc}{\newcommand}

\nc{\Id}{\text{Id}}
\nc{\la}{\lambda}

\begin{document}


\centerline{\bf LORENTZ GROUPS OF CYCLOTOMIC EXTENSIONS}

\bigskip\bigskip

\centerline{\it Relativity and Reciprocity}

\  

\centerline{Vadim Schechtman}

\

\centerline{March 16, 2021}

\bigskip\bigskip

\centerline{\bf Abstract}

\

In this (mostly historical) note 
we show how a unified Kummer-Artin-Schreier sequence from [W], [SOS] may be recovered from 
the relativistic velocity addition law. 


\

\centerline{\bf Introduction}

\

In the beginning of XX-th century, almost at the same time, two articles have appeared. 

One of them has been published in {\it Annalen der Physik}, and became very celebrated. 

The other one appeared in {\it Mathematischen Annalen}; it received the 
prize from the {\it K\"oniglische Gesellschaft der Wissenschaften zu G\"ottingen}. 
It was written by a number theorist Philipp Friedrich Furtw\"angler\footnote{{\it "der zu den bedeutensten Zahlentheoretiker seiner Zeit geh\"orte"}} 
(see [F], [H]) and  
contained a proof of a reciprocity law for $l$-th powers envisioned by Hilbert, generalizing the classical quadratic reciprocity. 

Later on (see [SS]) it was remarked that the Furtw\"angler's definitions contain implicitly certain 
group scheme which approximates between the multiplicative group  $\mu_p$ of $p$-th roots of unity (the kernel of the Kummer map) and its additive analogue 
$\alpha_p$ (the kernel of the Artin-Schreier map)
(see also [ST], Exercises 1 and 2). 

The first mentioned paper is [E], where the author discusses the consequences of the postulate 
that the speed of light $c$ remains constant in moving frames.


It is based on what became known as a relativistic formula for velocity addition. 

In the present note we remark  that a unified Kummer-Artin-Schreier group scheme may be 
extracted from this addition formula, so as the case $c$ finite ("Lorentzian") corresponds to Kummer, 
whereas $c = \infty$ ("Galilean") corresponds to Artin-Schreier.   

The reader may notice  a connection between the relativistic addition formula 
and symmetric functions described in 1.1.


\

{\bf Acknowledgement.} This note was inspired by [B].  
I am grateful to Y.Bazaliy, J.Tapia and B.Toen for consultations.

\vspace{1cm}

\centerline{\bf \S 1. Relativistic velocity addition law}

\bigskip\bigskip

At the moment we will work over a ground commutative ring $k$; in 1.3 below 
we will suppose  that $k = \BR$.

{\bf 1.1. Velocity addition and symmetric functions.} Denote
$$
u\oplus v = u\oplus_h v = \frac{u + v}{1 + h^2uv} = \frac{\sigma_1(u, v)}{1 + h^2\sigma_2(u, v)}
\eqno{(1.1.1)}
$$
This operation is commutative (obvious) and associative:
$$
(u\oplus v) \oplus w = u\oplus (v \oplus w) = 
\frac{\sigma_1(u,v,w) + h^2\sigma_3(u,v,w)}{1 + h^4\sigma_2(u,v,w)}
\eqno{(1.1.2)}
$$
More generally,
$$
u\oplus v \oplus w \oplus t = \frac{\sigma_1(u,v,w,t) + h^2\sigma_3(u,v,w,t)}{1 + h^2\sigma_2(u,v,w,t) + h^4\sigma_4(u,v,w,t)},
\eqno{(1.1.3)}
$$
etc.

The operation $\oplus_h$ has the neutral element $0$, and the inverse $u\mapsto -u$.

We get a rational group law to be denoted by $\BG_h$. Thus $\BG_h = \Spec k[x]$, and the 
group law is a rational map
$$
\oplus_h: \ \BG_h\times \BG_h --\lra \BG_h.
$$

{\bf 1.2. Pseudo-orthogonal $2\times 2$-matrices.}
Let
$$
A(u) = A_h(u) = \left(\begin{matrix} 1 & - h^2u\\
- u & 1\end{matrix}\right).
$$
Note that $A_0(u)$ is just a lower triangular matrix. 

Introduce a scalar product $( , )_h$ on $V = k^2$ by 
$$
((a,b), (a',b'))_h = h^2aa' - bb'.
\eqno{(1.2.1)}
$$
Then the lines (or columns) of $A(u)$ are orthogonal with respect to this scalar product, we can write
$$
A(u)\in O(2, ( , )_h).
$$
We have 
$$
A(u)A(v) = (1 + h^2uv)A(u\oplus_h v).
\eqno{(1.2.2)}
$$
Taking determinants, we get
$$
a(u)a(v) = ( 1 + h^2uv) a(u\oplus v)
$$
where
$$
a(u) = a_h(u) = \det A(u) = 1 - h^2u^2.
$$
In other words, the function
$$
b(u, v) = b_h(u, v) = 1 + h^2uv = \sigma_0(u,v) + h^2\sigma_2(u,v)
$$
is a coboundary: 
$$
\frac{a(u\oplus v)}{a(u)a(v)} = \frac{1}{b(u,v)^2}.
\eqno{(1.2.3)}
$$
It follows that if we define
$$
B_h(u) = \frac{A_h(u)}{a_h(u)^{1/2}}
$$
then 
$$
B_h(u)\in SO(2, (,)_h) = \BL_h \subset SL_2
$$
and
$$
B_h(u)B_h(v) = B_h(u\oplus_h v).
\eqno{(1.2.4)}
$$
For $h = 0$ the matrices $B_0(u) = A_0(u)$ generate a group 
$\BL_0$ of lower unipotent matrices isomorphic to $\BG_a$. 

\

{\it Physical interpretation}

\

Consider the relativistic $1|1$-dimensional space-time where the speed of light is $c = 1/h$. 
Denote
$$
L_c(u) := B_h(u).
$$
If $(x, t)$ are coordinates of an event in a frame at rest then
$$
\left(\begin{matrix} x'\\ t'
\end{matrix}\right) = L_c(u)\left(\begin{matrix} x\\ t
\end{matrix}\right)
$$
are its coordinates in the frame moving at the velocity $u$, cf. for example [B].  

\

{\bf 1.3. Positivity.} Let $k = \BR$. Suppose that $h > 0$, and let 
$$
h = \frac{1}{c}, \ c > 0
$$
(so $c$ is the "velocity of light"). Then $\oplus_h = \oplus_c$ gives rise to a group law on 
$$
I_c := (-c, c),
$$
i.e. $I_c$ becomes a Lie group.

The arrow $\beta_h$ (see (2.1.2) below) is an isomorphism of real Lie groups 
$$
I_{c} \iso \BR_{>0}^*
$$
Thus we have got a family of embeddings of  Lie groups
$$
\iota(c): \BR_{>0}^* \iso I_c \hra SL_2(\BR) 
$$
which degenerates in the limit $c\lra\infty$ to an oryspheric (lower triangular) subgroup
$$
\BR \isom U_-\hra SL_2(\BR).
$$
Replacing $c$ by $ic$ we get a family of compact tori 
$$
S^1\isom T_c\hra SL_2(\BR).
$$

\

\centerline{\bf \S 2. Relation to $\BG_m$}

\

{\bf 2.1. } We claim that for $h\in k^*$ the group law $\oplus_h$ 
is (essentially) isomorphic to $\BG_m$. 


We have
$$
1 + h\frac{u  + v}{1 + h^2uv} = \frac{(1 + hu)(1 + hv)}{1 + h^2uv},
$$
therefore
$$
\frac{1 + h(u\oplus_{h}v)}{\sqrt{a_{h}(u\oplus_{h}v)}} = 
\frac{1 + hu}{\sqrt{a_{h}(u)}}\cdot \frac{1 + hv}{\sqrt{a_{h}(v)}}
\eqno{(2.1.1)}
$$
This means that if we define
$$
\beta_h:\ \BL_{h} \lra \BG_m
$$
by 
$$
\beta_h(u) = \frac{1 + hu}{\sqrt{1 - h^2u^2}} = \sqrt{\frac{1 + hu}{1 - hu}}
\eqno{(2.1.2)}
$$ 
(on the level of points) then $\beta_h$ is a morphism of group laws. 


The inverse map 
$$
\beta_h^{-1}:\ \BG_m \lra \BL_h
$$
is
$$
v\mapsto \frac{1}{h}\frac{v^2 - 1}{v^2 + 1},
\eqno{(2.1.3)}
$$
well defined if $h\in k^*$. 

In fact we can get the Lorentz group law (1.1.1) by transfering the group law from $\BG_m$  
using the maps $\phi_h, \phi_{h}^{-1}$.

\

{\it  Relation to tori}

\

The other way around, we can start from matrices $A_h(u)$.

We remark that $A_h(u)$ has eigenvalues 
$\lambda_\pm(u) = 1 \pm hu$ with eigenvectors 
$$
v_\pm = \left(\begin{matrix}\mp h\\ 1
\end{matrix}\right).
$$
This implies that 
$$
C_h^{-1}A_h(u)C_h = \left(\begin{matrix} 1 + hu & 0 \\ 0 & 1 - hu 
\end{matrix}\right)
$$
where
$$
C_h = \left(\begin{matrix} - h  & h \\ 1 & 1  
\end{matrix}\right).
$$
Therefore the subgroup $\BL_h\subset SL_2$ is conjugated be means of $C_h$ to the standard 
maximal torus, and the map $\beta_h$ is the conjugation.

Thus we've got a family of tori which degenerates as $h\lra 0$ to an oryspheric subgroup 
$\BL_0 = U_-$.

\

{\bf 2.2. SOS group.} 
Cf. [W], [SOS], [SS]. 

Let 
$$
\BS_h = \BG'_h = \Spec k[x, (1 + hx)^{-1}]
$$
We have a morphism
$$
\alpha_h:\ \BS_h\lra \BG_m = \Spec k[z, z^{-1}],\ z\mapsto 1 + hx  
$$
which is an isomorphism of $h\in k^*$. On the level of points:
$$
\alpha_h(u) = 1 + hu,
$$
$$
\alpha^{-1}_h(v) = \frac{1}{h}(v - 1).
$$ 

Using it we can transfer the group law on $\BG_m$ to that on $\BS_h$. Since
$$
(1 + hx)(1 + hy) = 1 + h(x + y + hxy),
$$
the resulting group law on $\BS_h$ will be (on the level of points) 
$$
u \oplus'_h v = u + v + huv
\eqno{(2.2.1)}
$$
which makes sense for any $h\in k$. The inverse:
$$ 
u\mapsto - \frac{u}{1 + hu}.
$$

\

{\it The Kummer map}

\

For any $n\in \BZ$ we have a group map
$$
n:\ \BG_m \lra \BG_m,\ v\mapsto v^n.
$$
We define a group map
$$
\psi_n:\ \BS_h \lra \BS_{h^n}
$$
by transferring, so that the square 
$$
\begin{matrix} 
\BS_h & \overset{\psi_n}\lra & \BS_{h^n}\\
\alpha_h\downarrow & & \downarrow\alpha_{h^n}\\
\BG_m & \overset{n}\lra & \BG_m\\
\end{matrix}
$$
would be commutative, whence the formula
$$
\psi_n(u) = \frac{(hu + 1)^n - 1}{h^n},
\eqno{(2.2.2)}
$$ 
cf. [SOS]. 

\

{\bf 2.3. The Kummer map for $\BL_h$.} We play the same game with $\BL_h$. 

Namely a morphism of groups $\phi_n$ included into a commutative diagram
$$
\begin{matrix} 
\BL_h & \overset{\phi_n}\lra & \BL_{h^n}\\
\beta_h\downarrow & & \downarrow\beta_{h^n}\\
\BG_m & \overset{n}\lra & \BG_m\\
\end{matrix}
$$
is defined by
$$
\phi_n(u) = \frac{1}{h^n}\frac{(1 + nh)^n - (1 - nh)^n}{(1 + nh)^n + (1 - nh)^n}.
$$

\

\newpage

\centerline{\bf \S 3. Roots of unity and Artin-Schreier}

\

\

In 3.1, 3.2 below we recall [SOS]; in 3.3 we describe its analogue for the Lorentz group.

\

{\bf 3.1. $p$-cyclotomic extension of $\BQ_p$. } Together with [SOS] consider a ring $k = \BZ_p[\zeta_p]$, 
$\zeta = \zeta_p = e^{2\pi i/p}$, $p$ being a prime number. Thus $k$ is the ring of integers 
in $\BQ_p(\zeta_p)$. 

Let $h = \zeta - 1 \in k$. 

We have a totally ramified extension of degree $p - 1$,
$$
\BZ_p \subset k = \BZ_p[y]/(1 + y + \ldots + y^{p-1}), 
$$
$k$ is a discrete valuation ring with a uniformising parameter $h$ and quotient field $k/(h) = \BF_p$, cf. [Se], Ch. IV, Prop. 17.

We have
$$
h^{p-1} = w p
\eqno{(3.1.1)}
$$
where $w\in k^*$ and
$$
w \equiv - 1 \mod h.
\eqno{(3.1.2)}
$$

{\bf 3.1.1. Example.} Let $p  = 3$; then:
$$
\zeta^2 + \zeta + 1 = 0,
$$
$$
h^2 = (\zeta - 1)^2 = - 3\zeta,
$$
$$
w = - \zeta = (1 + \zeta)^{-1} = - 1 - h.
$$

\

(a) {\it The group $\BS_h$}

\

{\bf 3.2.} Consider the map 
$$
\psi_p:\ \BS_h\lra \BS_{h^p},
$$ 
$$
\psi_p(u) = \frac{(hu + 1)^p - 1}{h^p}
$$
Due to $(3.1.1)$ this map is defined over $k$.

On the other hand due to $(3.1.2)$ its special fiber
$$
\psi_{p,s} := \psi_p\otimes_k\BF_p: \BG_{a,\BF_p} \lra \BG_{a,\BF_p}
$$
coincides with the Artin-Schreier map
$$
\psi_{p,s}(u) = \wp(u) := u^p - u,
$$
cf. [SOS], (2.2.2), (2.2.3).

\

(b) {\it The group $\BL_h$}

\
 
{\bf 3.3.} Similarly consider the map 
$$
\phi_p:\ \BL_h\lra \BL_{h^p},
$$ 
$$
\phi_p(u) = \frac{1}{h^p}\frac{(1 + ph)^n - (1 - ph)^n}{(1 + ph)^n + (1 - ph)^n}.
$$
Again due to $(3.1.2)$ its special fiber at $h = 0$ coincides with the Artin - Schreier map. 
 

\ 


\centerline{\bf \S 4. Remarks on $q$-deformed tori}

\

In this Section we use our viewpoint to clarify somewhat the contents of [MRT], 6.4.1.

{\bf 4.1. Filtered circle.}   In [MRT], 6.4.1 the authors introduce an object 
$\bfS^1_{\Fil,\BZ}$ 
which they call the ''filtered circle over $\BZ$''. They start with the formal group law $\BS_h$ (2.2.1) over $\BA^1_\BZ$  which they denote $\BG$, $h$ (their $\lambda$) being the coordinate along $\BA^1$. 

Let $\CR = \Dist(\BG)$ denote the algebra of distributons - the dual to the algebra of functions 
on $\BG$; it is a filtered commutative and cocommutative Hopf algebra. They 
define a filtered group scheme 
$$
\bfH_\BZ = \Spec(Rees(\CR))
$$
where $Rees$ denotes the Rees construction. By definition $\bfS^1_{\Fil,\BZ}$ is its classifying stack 
$$
\bfS^1_{\Fil,\BZ} = B\bfH_\BZ.
$$
This construction can be $q$-deformed. To do this, the authors remark that the fiber at $h = 1$, 
$\CR_1 = \Dist(\BG_1)$ is isomorphic to the algebra  of polynomials $f(x)\in\BQ[x]$ such that 
$f(\BZ)\subset \BZ$ (I am grateful to B.Toen who explained to me the above statement). The last algebra admits a $q$-deformation $\CR_{1,q}$ introduced in [HH].
  
This allows to define a $q$-deformation $\CR_q$ of $\CR$. 


It was remarked in [HH] that the algebra $\CR_{1,q}$  is isomorphic to the standard  Cartan 
subalgebra $K_q$ of the Lusztig's divided power quantum group $U_q\fsl_2$. 
 
\

{\bf 4.2. Lorentzian counterpart.}  
We may replace in the above discussion $\BS_h$ by $\BL_h$ which by definition lies inside $SL_2$; 
so we consider a model $SL_{2,\BZ}$ of $SL_2$ over $\BZ$, and $\{\BL_h\}$ is a family of conjugated tori inside it.

We can consider a $q$-deformation of $\Dist(SL_{2,\BZ})$, i.e. the corresponding quantum group $U_q\fsl_2$, 
and it would be natural to think that a family of quantum tori from 4.1 may be identified with  
a $q$-deformation $\BL_{h,q} \subset U_q\fsl_2$. However a realization of this idea meets 
some difficulties at the moment.


\

{\bf 4.3. Towards a filtered $S^3$.} The Lie group $SL_2(\BC)$ is homotopy equivalent to $SU(2)$ (the Weyl's unitary trick) which is in turn homeomorphic to the sphere $S^3$. 

This rises a question whether the whole packet of [MRT] and [R] admits a noncommutative ''quaternionic'' version, with $S^1$ replaced by  $S^3$, cf. [Sch].

\

\centerline{\bf Literature}

\bigskip\bigskip

[E] A.Einstein, Zur Elektrodynamik bewegter K\"orper, {\it Ann. der Physik} {\bf  17} (1905), 891 - 921.

[F] Ph.Furtw\"angler, \"Uber die Reziprot\"atsgesetze zwischen $l$-ten Potenzresten in algebr\"aischen 
Zahlk\"orpers, {\it Abh. K\"oniglischen Ges. Wiss. G\"ottingen} {\bf 2} (1902), Heft 3, 1 - 82; 
{\it Math. Annalen} {\bf 58} (1904), 1-50.

[B] Yaroslav Bazaliy, Three zoom talks on Special relativity, Toulouse "Oxford" seminar, Fall 2020.

[GMS] I.M.Gelfand, R.A.Minlos, Z.Ya.Shapiro, Representations of the rotation and Lorentz groups and their applications

[HH] Nate Harman, Sam Hopkins, Quantum integer-valued polynomials, arXiv:1601.06110. 

[H] Helmut Hasse, History of Class field thery, in: Algebraic Number Theory, J.W.S.Cassels, A.Fr\"ohlich 
(eds). 1967, 266-279.

[MRT] Tasos Moulinos, Marco Robalo, Bertrand Toen, A universal HKR 
\newline theorem, arXiv:1906.0011

[M] N.David Mermin, It's about time

[R] Arpon Raksit, Hochschild homology and the derived de Rham complex revisited, arXiv:2007.02576. 

[Sch] V.Schechtman, Moufang loops and toric surfaces, arXiv:2102.12868.

[SOS] T.Sekiguchi, F.Oort, N.Suwa, On the deformation of Artin-Schreier to Kummer, {\it Ann. Sci. ENS} {\bf 22} (1989), 345 - 375.

[Se] J.-P.Serre, Corps locaux

[ST] J.-P.Serre, J.Tate, Exercises, in: Algebraic Number Theory, J.W.S.Cassels, A.Fr\"ohlich 
(eds). 1967, 266-279.

[SS] Noriyuki Suwa, Tsutomu Sekiguchi, Th\'eorie de Kummer-Artin-Schreier et applications, 
{\it J. Th. Nombr. de Bordeaux} {\bf 7} (1995), 177-189.

[W] William C.Waterhouse, A unified Kummer-Artin-Schreier sequence, {\it Math. Ann. } {\bf 277} (1987), 
447 - 451.

\

\end{document}